\documentclass[twocolumn]{autart}    % Enable this line and disable the 
\usepackage{cite}  
\PassOptionsToPackage{pdfpagelabels=false}{hyperref} 
\usepackage{graphicx}
\usepackage{amsmath}
\usepackage{mathtools, amssymb}
\usepackage{varioref}
\labelformat{equation}{(#1)}
\usepackage{amsfonts}

\usepackage{enumerate}
\usepackage{color}
\usepackage{mysymbol}

\usepackage{empheq}

\usepackage{epstopdf}
\usepackage[dvipsnames]{xcolor}
 \usepackage{tikz}%,pgfplots}
\usepackage{color} % change text color        
\usetikzlibrary{matrix} % for block alignment
\usetikzlibrary{arrows} % for arrow heads
\usetikzlibrary{calc} % for manimulation of coordinates
\usetikzlibrary{shapes}
\usetikzlibrary{arrows}
\usetikzlibrary{positioning}
\usetikzlibrary{backgrounds, fit}

\usepackage{url}

\theoremstyle{plain}
\newtheorem{lemma}{Lemma}
\newtheorem{proposition}{Proposition}

\theoremstyle{definition}
\newtheorem{asmptn}{Assumption}
\newtheorem{remark}{Remark}

\newtheorem{example}{Example}

\usepackage{algorithm}
\usepackage{algpseudocode}

\usepackage{standalone}

\usepackage{enumerate}

\newcommand{\figref}[1]{Fig.~\ref{fig:#1}}

\newcommand{\qedsymbol}{\hspace*{\fill}~$\blacksquare$}
\newcommand{\qedopensymbol}{\hspace*{\fill}~$\square$}

\renewcommand{\epsilon}{\varepsilon}

\begin{document}

\begin{frontmatter}

\title{
Random Access Design for Wireless Control Systems
}
\thanks{This work was supported in part by NSF CNS-0931239, and by TerraSwarm, one of six centers of STARnet, a Semiconductor Research Corporation program sponsored by MARCO and DARPA.}
%
%\thanksref{footnoteinfo}
%} % Title, preferably not more 
																								
\author{Konstantinos Gatsis}\ead{kgatsis@seas.upenn.edu},% Add the 
\author{Alejandro Ribeiro}\ead{aribeiro@seas.upenn.edu},% e-mail address 
\author{George J. Pappas}\ead{pappasg@seas.upenn.edu}% (ead) as shown

\address{Department of Electrical and Systems Engineering, University of Pennsylvania, 200 South 33rd Street, Philadelphia, PA 19104}  % Please supply                                              
%\address[Rome]{Senate House, Rome}             % full addresses
%\address[Baiae]{The White House, Baiae}        % here.
          
                             % keyword list or with the 
                                          % help of the Automatica 
                                          % keyword wizard

\begin{abstract}                          % Abstract of not more than 200 words.
Interferences arising between wireless sensor-actuator systems communicating over shared wireless channels adversely affect closed loop control performance. To mitigate this problem we design appropriate channel access policies for wireless control systems subject to channel fading. The design is posed as an optimization problem where the total transmit power of the sensors is minimized while desired control performance is guaranteed for each involved control loop. Control performance is abstracted as a desired expected decrease rate of a given Lyapunov function for each loop. We prove that the optimal channel access policies are decoupled and, intuitively, each sensor balances the gains from transmitting to its actuator with the negative interference effect on all other control loops. Moreover the optimal policies are of a threshold nature, that is, a sensor transmits only under favorable local fading conditions. Finally, we show that the optimal policies can be computed by a distributed iterative procedure which does not require coordination between the sensors.
%
%We consider multiple sensors randomly accessing a shared wireless channel to transmit measurements of their respective plants to a controller. The channel access policies of the sensors are appropriately designed to mitigate packet collisions from simultaneous transmissions and to exploit transmission opportunities due to random channel fading conditions. This design aims to minimize the total transmit power of the sensors while control performance for all control loops is guaranteed. Control performance is abstracted as a desired expected decrease rate of a given Lyapunov function for each loop. By Lagrange duality arguments, the optimal access policies are shown to be decoupled among sensors and of a threshold form, that is, each sensor transmits only if its local channel state is favorable enough. Moreover, to compute the optimal policies we develop an iterative procedure which is distributed and easily implementable as it does not require coordination between the sensors.
%
\end{abstract}

\end{frontmatter}

\section{Introduction}

Wireless sensors are ubiquitous in modern smart infrastructures where they are deployed to monitor and control physical processes in our homes, urban environments, and industrial plants. This abundance of wireless devices however also creates an increase in the wireless interferences arising between transmissions over the shared wireless medium. The development of decentralized communication mechanisms that can mitigate these interference effects and guarantee closed loop control performance arises as an important research direction.

%creates the need to efficiently manage the available wireless medium to be shared among these devices in way that control performance guarantees are met.

%The abundance of wireless sensing devices in modern control applications, for example, smart homes, industrial environments, and urban infrastructure, creates a need for sharing the available wireless medium in an efficient manner. The primary goal of such channel access mechanisms is to provide control performance guarantees. Besides this goal, the development of mechanisms that are also easily implementable, e.g., decentralized, arises as an important research direction.

The prevalent approach to the problem sharing a wireless communication medium in networked control systems is centralized scheduling which guarantees no interferences. Static scheduling for example specifies that sensors transmit in some predefined periodically repeating sequence such as round-robin and this sequence is designed to meet control objectives, see, e.g.,~\cite{Branicky_stability, Hristu_shared_feedback, LeNy_resource_LQR}. Deriving optimal scheduling sequences is recognized as a hard combinatorial problem~\cite{Scheduling_control_combinatorics, Gupta_Hassibi_sensor_selection}. Scheduling can also be dynamic, where at each time step a central network coordination authority decides which device gets access to the medium. This dynamic decision may be stochastic~\cite{Gupta_Hassibi_sensor_selection}, based on plant state information~\cite{Walsh_stability, Donkers_switched}, or based on the wireless channel conditions~\cite{GatsisEtal15}.

In contrast to centralized scheduling, decentralized mechanisms where sensors independently decide access to the shared wireless medium are easier to implement. They do not require predesigned sequences of how sensors access the medium, or a central authority to take scheduling decisions. The drawback of this decentralized approach however is that packet collisions can occur from simultaneously transmitting sensors, resulting in lost packets and control performance degradation. Hence sensor access policies need to be appropriately designed to mitigate these effects. We consider specifically a random access mechanism where each sensor independently and randomly decides whether to transmit plant state measurements over a shared channel to an access point/controller (\figref{random_access}). 

Control under random access communication mechanisms has drawn limited attention, to the best of our knowledge. Comparisons between different medium access mechanisms for networked control systems and the impact of packet collisions in stability and control performance have been considered either in numerical simulations~\cite{Goldsmith_MAC_NCS, Ramesh_Johansson} or analytically in simple cases~\cite{Rabi_Johansson_Estimation_MAC, Aloha_event}. These include random access mechanisms and related Aloha-like schemes, where after a packet collision the involved sensors wait for a random time interval and retransmit. Stability conditions under packet collisions were examined in \cite{Nesic_io_stability, NCS_Aloha}. In contrast to these works, our goal is to design the medium access mechanism so that desired control performance is guaranteed. Besides closed loop control, optimal remote estimation over collision channels is considered recently in \cite{Martins_estimation_shared}.

We pose the design of channel access policies for multiple control loops over a shared wireless channel as an optimization problem (Section~\ref{sec:problem}). The goal is to satisfy a control performance requirement for each control loop while minimizing the total expected transmit power expenditures of the sensors. We propose a Lyapunov-like control performance abstraction, motivated from our work on centralized scheduling~\cite{GatsisEtal15}. Each control system is abstracted via a given Lyapunov function which is desired to decrease at a predefined rate and in expectation due to the random packet losses and collisions on the shared medium. These control requirements are then shown to be equivalent to a minimum packet success rate on each link. 

Besides mitigating packet collisions, sensors can exploit channel fading state information. Fading refers to large unpredictable variations in wireless channel transferences~\cite[Ch.~3,4]{Goldsmith_Wireless_Communications} which in our setup affect the likelihood of successful packet decoding at the receiver. % by a known complementary error function. 
This communication model has been used in estimation and control applications~\cite{GatsisEtal14, GatsisEtal15, Quevedo_Kalman} but not under a random access mechanism.
%\KG{I have to emphasize the importance of random wireless channel. What does it offer me eg in power adaptation or centralized scheduling?}
%
We design sensor access policies that adapt to channel states, allowing to, e.g., transmit at higher rates under channel conditions with higher packet success. 
%We examine the design of sensor access rates that satisfy the Lyapunov control performance requirements and minimize the average transmit power of the sensors. 
In preliminary work presented in~\cite{c_GatsisEtal15} we considered again random access wireless sensors but employing simpler policies, in particular policies that do not adapt to channel states online.

Based on Lagrange duality arguments we characterize the structure of the optimal sensor access policies (Section~\ref{sec:solution}). We show that the optimal policies are of a threshold nature, that is, each sensor transmits only when its corresponding channel state is favorable enough and avoids transmission otherwise. Moreover we reveal an intuitive decoupling of the policies among sensors. Each sensor should select its channel threshold in a way that balances the control performance of its own closed loop with the collective negative effect it has on all other control loops due to collisions. Decentralized policies with similar structures are also known to be optimal for general wireless random access communication networks~\cite{Adireddy_Tong, Qin_Berry, Hu_Ribeiro_RA}. The context differs however, since in these works the objective is thoughput-based utility functions in contrast to the packet success rates used for control performance here. 
%\cite{Al-Harthi_etal} has the same collision model with CSI !

In Section~\ref{sec:optimization} we derive an iterative procedure to compute the optimal access policies. The procedure is easy to implement in our architecture as it does not require the sensors to coordinate among themselves, or to know what control performance the other sensors try to achieve. Technically the procedure optimizes the Lagrange dual problem, and relies on the common access point to compute the optimal dual variables and provide them to the sensors via the reverse channel. We conclude with a numerical example and some remarks (Sections~\ref{sec:simulations},~\ref{sec:conclusion}).

\section{System Description}\label{sec:problem}

\begin{figure}[!t]
\centering
\resizebox{1.0\columnwidth}{!}{
\includegraphics[width = \columnwidth]{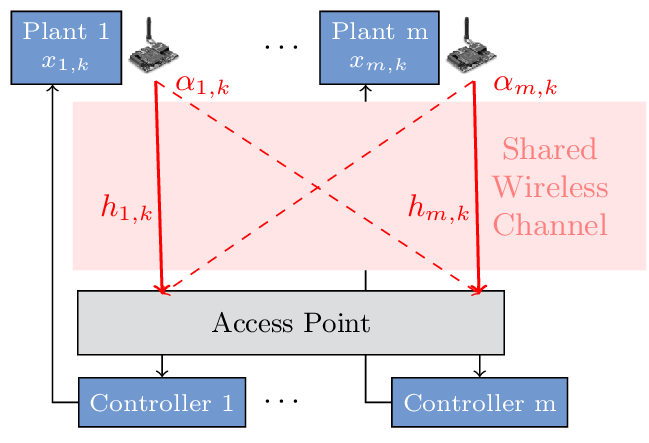}
}
\vspace{-0pt}
\caption{
Random access architecture for $m$ control loops over a shared wireless medium. Each sensor $i$ randomly transmits with probability $\alpha_{i,k}$ at time $k$ to a common access point computing the plant control inputs. If only sensor $i$ transmits, the successful decoding probability depends on local channel conditions $h_{i,k}$. If other sensors transmit at the same time a collision might occur at sensor $i$'s transmission, rendering $i$'s packet lost.% The goal is to design sensor access policies that guarantee control performance for all control loops.
}
\vspace{-0pt}
\label{fig:random_access}
\end{figure}

We consider a wireless control architecture where $m$ independent plants are controlled over a shared wireless medium. Each sensor $i$ ($i=1,2,...,m$) transmits measurements of plant $i$ to an access point responsible for computing the plant control inputs. Packet collisions might arise on the shared medium between simultaneously transmitting sensors. See \figref{random_access} for an illustration. We are interested in designing a mechanism for each sensor to independently decide whether to access the medium (random access) in a way that guarantees desirable control performance for all control systems.

Our goal is to design communication policies under the assumption that the dynamics for all $m$ control systems have been designed a priori and  are therefore fixed and independent of the communication policy. As a consequence, the system evolution is described by a switched model that depends on whether the controller manages to reach the access point or not. Thus, if we use $\gamma_{i,k} \in \{0,1\}$ to indicate the success of the transmission at time $k$ for link/system $i$ and assume the system is linear and time invariant, we can model its evolution by the switched system
\begin{eqnarray}\label{eq:system}
	x_{i,k+1} = \left\{ \begin{array}{ll}
	A_{c,i} \, x_{i,k} + w_{i,k}, &\text{ if } \gamma_{i,k}=1, \\
	A_{o,i} \, x_{i,k} + w_{i,k}, &\text{ if } \gamma_{i,k}=0.
	\end{array}\right. 
\end{eqnarray}
Here $x_{i,k} \in \reals^{n_i}$ denotes the state of control system $i$ at each time $k$, which can in general include both plant and controller states -- see, e.g., Example~\ref{ex:system}. At a successful transmission the system dynamics are described by the matrix $A_{c,i} \in \reals^{n_i \times n_i}$, where 'c' stands for closed-loop, and otherwise by $A_{o,i} \in \reals^{n_i \times n_i}$, where 'o' stands for open-loop. We assume that $A_{c,i} $ is asymptotically stable, implying that if system $i$ successfully transmits at each slot the state evolution of $x_{i,k}$ is stable. The open loop matrix $A_{o,i} $ may be unstable. The additive terms $w_{i,k}$ model an independent (both across time $k$ for each system $i$, and across systems) identically distributed (i.i.d.) noise process with mean zero and covariance $W_i\succeq 0$. An example of such a networked control system is presented next.

\begin{example}\label{ex:system}
Suppose each closed loop $i$ consists of a linear plant and a linear output of the form
\begin{eqnarray} \label{eq:plant}
	%\begin{split}
	&x_{i,k+1} &= A_i x_{i,k} + B_i u_{i,k} + w_{i,k}, \\
	&y_{i,k} &= C_i x_{i,k} + v_{i,k}, %\notag
	%\end{split}
\end{eqnarray}
where $\{w_{i,k}, \, k \geq 0\}$ and $\{v_{i,k}, \, k \geq 0\}$ are i.i.d. Gaussian disturbance and measurement noise respectively. Each wireless sensor $i$ transmits the output measurement $y_{i,k}$ to the controller. A dynamic control law adapted to the packet drops keeps a local controller state $z_{i,k}$, 
\begin{eqnarray} \label{eq:dynamic_controller_1}
	%\begin{split}
	&z_{i,k+1} &= F_i \, z_{i,k} + \gamma_{i,k} \, ( F_{c,i} \, z_{i,k} + G_i \, y_{i,k})
	%\end{split}
\end{eqnarray}
which may for example represent a local estimate of the plant state~\cite{Hespanha_survey}, and applies plant input $u_{i,k}$ as
\begin{eqnarray} \label{eq:dynamic_controller_2}
	%\begin{split}
	&u_{i,k} &= K_i \, z_{i,k} + \gamma_{i,k} \, ( K_{c,i} \, z_{i,k} + L_i \, y_{i,k}). %\notag
	%\end{split}
\end{eqnarray}
In other words, the controller updates appropriately the local state and input whenever a measurement is received. The overall closed loop system is obtained by joining plant and controller states into
\begin{eqnarray} \label{eq:overall_system}
\left[ \begin{array}{c}
	x_{i,k+1} \\ z_{i,k+1} \end{array}
	\right] = 
	\left[ \begin{array}{cc} 
	A_i+ \gamma_{i,k} B_i L_i C_i &B_i K_i+ \gamma_{i,k}B_i K_{c,i} \\ \gamma_{i,k}G_i C_i &F_i+ \gamma_{i,k} F_{c,i}
	\end{array} \right] \notag\\
	\qquad \quad \quad \; \;\cdot \left[ \begin{array}{c}
	x_{i,k} \\ z_{i,k} \end{array}
	\right] + 
	\left[ \begin{array}{cc} 
	I &\gamma_{i,k} B_i L_i C_i \\ 0 &\gamma_{i,k} G_i  
	\end{array} \right]
	\left[ \begin{array}{c}
	w_{i,k} \\ v_{i,k} \end{array}
	\right] 
\end{eqnarray}
%
%
%\begin{eqnarray} \label{eq:overall_system}
%%
%\left[ \begin{array}{c}
	%x_{i,k+1} \\ z_{i,k+1} \end{array}
	%\right] = 
	%%
	%\left[ \begin{array}{cc} 
	%A_i &B_i K_i \\ 0 &F_i
	%\end{array} \right] 
	%%
	%\left[ \begin{array}{c}
	%x_{i,k} \\ z_{i,k} \end{array}
	%\right] + 
	%\left[ \begin{array}{cc} 
	%I &0 \\ 0 &0 
	%\end{array} \right]
	%\left[ \begin{array}{c}
	%w_{i,k} \\ v_{i,k} \end{array}
	%\right] 
	%\notag\\
	%%
	%\left[ \begin{array}{c}
	%x_{i,k+1} \\ z_{i,k+1} \end{array}
	%\right] =  \left[ \begin{array}{cc} 
	%A_i+B_i L_i C_i &B_i (K_i+ K_{c,i}) \\ G_i C_i &F_i+ F_{c,i}
	%\end{array} \right] 
	%%
	%%
	%\left[ \begin{array}{c}
	%x_{i,k} \\ z_{i,k} \end{array}
	%\right] + 
	%\left[ \begin{array}{cc} 
	%I & B_i L_i C_i \\ 0 & G_i  
	%\end{array} \right]
	%\left[ \begin{array}{c}
	%w_{i,k} \\ v_{i,k} \end{array}
	%\right] 
%\end{eqnarray}
%%
which is of the form \eqref{eq:system}.
\qedopensymbol
\end{example}

%\subsection{Communication model}

The transmission success indicator variables $\gamma_{i,k}$ are random with a distribution that depends on the communication policy which here is supposed to be a slotted random access policy. Specifically, communication takes place in time slots generically indexed by $k$. At every slot $k$ each sensor $i$ transmits over the shared channel with some probability $\alpha_{i,k} \in [0,1]$ to be designed. A sensor's transmission might fail due to two reasons, packet decoding errors and packet collisions. A collision might be experienced on link $i$, thereby rendering packet $i$ lost, if some other sensor $j\neq i$ transmits in the same time slot. We assume that such a collision event occurs with constant probability $q_{ji} \in [0,1]$, {\it given} that both sensors $i, j$ transmit in the slot. Thus, the probability that sensor $i$'s transmission is free of collisions, i.e., that no other sensor transmits and causes collisions on link $i$, equals $\prod_{j\neq i} [ 1 - \alpha_{j,k} \, q_{ji} ]$. See Remark~\ref{rem:collision} for details of this collision model.

If sensor $i$ transmits and has a collision free time slot, the success of decoding the message at the access point/receiver depends on the randomly varying channel conditions on link $i$. Denote then by $h_{i,k} \in \reals_+$ the current channel fading conditions for link $i$ at time $k$. We adopt a block fading model~\cite[Ch.~4]{Goldsmith_Wireless_Communications} whereby channel states $\{ h_{i,k},\, k \geq 0 \}$ are assumed constant during each transmission slot $k$, but i.i.d. across time with distribution $\phi_i$. We also assume channel states are independent among systems $i$, a common assumption in the literature~\cite{Adireddy_Tong, Qin_Berry, Hu_Ribeiro_RA}, %Yu_Giannakis,
as well as independent of the plant process noise $w_{i,k}$. We let $q(h_{i,k})$ denote the probability of successful transmission given the current channel state $h_{i,k}$. For more details on this communication model see Remark~\ref{rem:fading}. The function $q: \reals_+ \rightarrow [0,1]$ is assumed to be continuous and strictly increasing, i.e., higher channel fading states imply higher packet success probability. %To sum up, if no other sensors transmit, system $i$ closes its loop at time step $k$ with probability $\alpha_{i,k} \, q(h_{i,k})$. 

Combining the effects of collisions and packet losses due to fading, the probability that a packet is successfully decoded at the access point can be written as
\begin{eqnarray}\label{eq:gamma_probability_conditioned}
	\mathbb{P}( \gamma_{i,k} = 1) = \alpha_{i,k} \, q(h_{i,k}) \; \prod_{j\neq i} \Big[ 1 - \alpha_{j,k} \, q_{ji} \Big].
\end{eqnarray}
This expression states that the probability of system $i$ in \eqref{eq:system} closing the loop at time $k$ equals the probability that transmission $i$ is successfully decoded at the receiver, multiplied by the probability that no other sensor $j \neq i$ is causing collisions on $i$th transmission. 
%The product**** in this expression is a consequence of the fact that all sensors independently decide to access the channel
%\footnote{*****Precisely, if the vector $\delta \in \{0,1\}^m$ indicates which of the $m$ sensors transmit at time $k$, we have $\mathbb{P}( \gamma_{i,k} =1 \given \delta) = \delta_{i} q(h_{i,k}) \prod_{j\neq i} [ 1 - q_{ji}]^{\delta_{j}}$. Taking expectation over the independent Bernoulli variables $\delta_i \sim \text{Bern}(\alpha_{i,k})$ yields \eqref{eq:gamma_probability}.}. 

Channel states reveal information about how easy it is for each sensor to successfully communicate, assuming no other sensor transmits. We assume that before deciding whether to transmit each sensor has access to its respective channel state and may adapt accordingly. For example sensor $i$ may transmit with higher or lower rate $\alpha_{i,k}$ under favorable or unfavorable channel states $h_{i,k}$ respectively. Hence we design policies that are measurable functions of the form $\alpha_{i,k} = \alpha_i (h_{i,k})$. Since channel states are i.i.d. over time we restrict attention to stationary policies, and drop the time index when not necessary. The set of all access policies for sensor $i$ is then
\begin{eqnarray}\label{eq:alpha_set_def}
	\mathcal{A}_i = \{ \alpha_i : \reals_+ \rightarrow [0,1] \}
\end{eqnarray}
and the vector $\alpha(.)$ of access policies for all sensors belongs in the Cartesian product space $\mathcal{A} = \mathcal{A}_1 \times \ldots \mathcal{A}_{m}$. For fixed sensor access policies, the probability of successful transmission on link $i$ can be expressed as
\begin{eqnarray}\label{eq:gamma_probability}
	\mathbb{P}( \gamma_{i,k} = 1) = \mathbb{E}_{h_i} [ \alpha_i(h_i) \, q(h_i)]\, \prod_{j\neq i} \Big[ 1 - \mathbb{E}_{h_j} [\alpha_j(h_j)] \, q_{ji} \Big]. \notag \\ \qquad
\end{eqnarray}
This expression follows from \eqref{eq:gamma_probability_conditioned} by taking expectation with respect to the channel states and using the independence of channels among systems. The expectation is well-defined as both functions $\alpha(.)$, $q(.)$ are measurable and bounded in $[0,1]$ hence integrable.%In \eqref{eq:gamma_probability} the communication parameters $q_{ji},\, i,j \in\{1,\ldots,m\}$ are given, and the variables to be designed are the sensor access policies $\alpha$.

We make the following technical assumption on the probability distribution of channel states, which holds true for practically considered models~\cite[Ch. 3]{Goldsmith_Wireless_Communications}.

\begin{asmptn}\label{as:non_atomic}
The distributions $\phi_i$ of channel states $\{ h_{i,k},\, k\geq 0 \}$ for all $i=1,\ldots,m$ are non-atomic, i.e., have a continuous distribution function on $\reals_+$.
\end{asmptn}

%The following technical assumption on the form of the function $q(h)$ will be helpful in the subsequent sections, and is verified for cases of practical interest as shown in \figref{q_vs_SNR}. 
%
%\begin{asmptn}\label{as:q_function}
%%
%The function $q: \reals_+ \rightarrow [0,1]$ is strictly increasing. 
%%satisfies:
%%%
%%\begin{enumerate}[(a)]
	%%\item $q(0) = 0$, (??)
	%%\item $q(h)$ is continuous, and strictly increasing when $q(h)>0$, i.e., for $h'>h$ it holds that $q(h')>q(h)>0$.
	%%%\item for any $\mu \in \reals_+$ and for almost all values $h \in \mathcal{H}$ the set $\argmin_{0 \leq p \leq p_{\max}} \, p - \mu \, q(h\,p)$ is a singleton.
%%\end{enumerate}
%%%
%\end{asmptn}

\begin{remark}\label{rem:collision}
Our collision model captured by the probabilities $q_{ji}$ 
%has been also recently considered in \cite{general_random_access} and 
subsumes: 
i) the conservative case where simultaneous transmissions certainly lead to collisions ($q_{ji} =1$) usually considered in control literature, e.g.,\cite{NCS_Aloha, Nesic_io_stability}, ii) the case where simultaneously transmitted packets are not always lost ($q_{ji} <1$), 
e.g., due to the capture phenomenon~\cite{Capture}, 
and iii) the asymmetric case where different sensors $j, \ell$ interfere differently on link $i$, e.g., due to their spatial configuration. %$0 < q_{ji} < q_{\ell i} \leq 1$ if sensor $\ell$ is closer than sensor $j$.
%For simplicity the values $q_{ji} \in (0,1]$ are assumed nonzero throughout our paper. The case where some sensors $j$ do not interfere with some link $i$, i.e., $q_{ji}=0$, is similarly handled. In that case the product in the packet success probability \eqref{eq:gamma_probability_conditioned} is only over the subset of sensors $j \neq i$ that affect link $i$.
%
\end{remark}

\begin{figure}[!t]
\centering
\includegraphics[width=0.9\columnwidth]{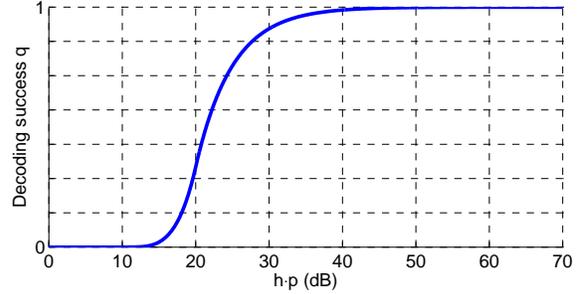}
\caption{Probability of successful decoding as a function of the received signal power level. Higher channel fading gains imply higher chance of packet success, hence more more favorable transmission opportunities.}
\label{fig:q_vs_SNR}
\end{figure}

\begin{remark}\label{rem:fading} 
The channel fading conditions $h_{i,k}$ change unpredictably over time~\cite[Ch.~3]{Goldsmith_Wireless_Communications} and affect the communication of the sensors. In particular if sensor $i$ transmits at a power level $p_i$, and assuming no other sensor transmits, the power level of the received signal is the product $h_{i,k} p_i$ of the current channel fading gain and transmit power. During high channel fading gains for sensor $i$ there is a higher received signal-to-noise ratio (SNR) at the access point/controller and consequently a higher chance to successfully decode the transmitted message. We let $q(h_{i,k})$ denote the packet success as a function of the channel state -- for more details on this model the reader is referred to~\cite{GatsisEtal14}. An illustration of this relationship is given in \figref{q_vs_SNR}. We suppose that each sensor $i$ has access to the channel state $h_{i,k}$ before deciding whether to transmit over the shared channel. For example this can be performed by short pilot signals sent from the access point to the sensors at the beginning of each time slot. 
\end{remark}

\section{Control Performance and Random Access Problem}

The random packet success on link $i$ modeled by \eqref{eq:gamma_probability} causes each control system $i$ in \eqref{eq:system} to switch in a random fashion between the two modes of operation (open and closed loop). As a result, the sensor access policies $\alpha(.)$ to be designed affect the performance of all control systems. The following proposition characterizes via a Lyapunov-like abstraction a connection between control performance and the packet success rate.

\begin{proposition}[Control performance abstraction]\label{prop:Lyapunov_abstraction}
Consider a switched linear system $i$ described by \eqref{eq:system} and $\gamma_{i,k}$ being a sequence of i.i.d. Bernoulli random variables, and a quadratic function $V_i(x_i) = x_i^T P_i x_i, \; x_i \in \reals^{n_i}$ with a positive definite matrix $P_i \succ 0$. 
%Suppose that the function $V_i(x_i)$ is a Lyapunov function for the closed loop mode of the system, i.e., $A_{c,i}^T \, P_i \, A_{c,i} \prec P_i$. 
Then the function $V_i(x_i)$ decreases with an expected rate $\rho_i<1$ at each step, i.e., we have
\begin{eqnarray}\label{eq:Lyapunov_condition}
	\mathbb{E} \left[ V_i(x_{i,k+1}) \given x_{i,k} \right]
	 \leq \rho_i \, V_i(x_{i,k}) + Tr(P_i W_i) 
\end{eqnarray}
for all $x_{i,k} \in \reals^{n_i}$, if and only if
\begin{eqnarray}\label{eq:gamma_requirement}
	\mathbb{P}(\gamma_{i,k}=1) \, \geq \, c_i, 
\end{eqnarray}
where $c_i \geq 0$ is computed by the semidefinite program
%%
%\begin{eqnarray}
	%c_i = \min\{ \theta \geq 0: \; \theta \, ( A_{o,i}^T P_i A_{o,i} - A_{c,i}^T P_i A_{c,i} ) \succeq ( A_{o,i}^T P_i A_{o,i} - \rho_i P_i )
%\end{eqnarray}
%
\begin{eqnarray}\label{eq:compute_c}
	c_i =  \min\{ &\theta & \geq 0:   \notag\\
	&\theta & A_{c,i}^T P_i A_{c,i} + (1-\theta) A_{o,i}^T P_i A_{o,i} \preceq \rho_i P_i \} .
\end{eqnarray}
\end{proposition}

\begin{pf*}{Proof.} %See Appendix~\ref{sec:Lyapunov_proof}.
The expectation over the next system state $x_{i,k+1}$ on the left hand side of \eqref{eq:Lyapunov_condition} accounts via \eqref{eq:system} for the randomness introduced by the process noise $w_{i,k}$ and the random success $\gamma_{i,k}$. In particular we have that
\begin{eqnarray}\label{eq:Lyapunov_evolution}
	&\mathbb{E} \left[ V_i(x_{i,k+1}) \given x_{i,k} \right] = 
	\mathbb{P}(\gamma_{i,k}=1) \; x_{i,k}^T A_{c,i}^T P_i A_{c,i} x_{i,k} \notag\\
	&+ \mathbb{P}(\gamma_{i,k}=0) \; x_{i,k}^T A_{o,i}^T P_i A_{o,i} x_{i,k} + Tr(P_i W_i).
\end{eqnarray}
Here we used the fact that the random variable $\gamma_{i,k}$ is independent of the system state $x_{i,k}$. Plugging \eqref{eq:Lyapunov_evolution} at the left hand side of \eqref{eq:Lyapunov_condition} we get 
\begin{eqnarray}\label{eq:Lyapunov_with_x}
	&&\mathbb{P}(\gamma_{i,k}=1) \; x_{i,k}^T A_{c,i}^T P_i A_{c,i} x_{i,k} \notag\\
	&&+ \mathbb{P}(\gamma_{i,k}=0) \; x_{i,k}^T A_{o,i}^T P_i A_{o,i} x_{i,k} \leq \, \rho_i x_{i,k}^T  P_i x_{i,k}.
\end{eqnarray}
%%
%for $x_{i,k} \neq 0$
%%
%\begin{eqnarray}\label{eq:Lyapunov_with_x}
	%\mathbb{P}(\gamma_{i,k}=1) \, \geq \, \frac{ x_{i,k}^T ( A_{o,i}^T P_i A_{o,i} - \rho_i P_i ) x_{i,k}  }
	%{ x_{i,k}^T ( A_{o,i}^T P_i A_{o,i} - A_{c,i}^T P_i A_{c,i} ) x_{i,k} }. 
%\end{eqnarray}
%%
Since condition \eqref{eq:Lyapunov_condition} needs to hold for all $x_{i,k} \in \reals^{n_i}$ we can rewrite \eqref{eq:Lyapunov_with_x} as a linear matrix inequality~\cite{Boyd_convex}
\begin{eqnarray}\label{eq:LMI}
	&&\mathbb{P}(\gamma_{i}=1)  A_{c,i}^T P_i A_{c,i} + (1 - \mathbb{P}(\gamma_{i}=1) )  A_{o,i}^T P_i A_{o,i} 
	\preceq \rho_i  P_i ,
	\notag\\ 
	&&\qquad
\end{eqnarray}
where we dropped the time indices from $\gamma_{i,k}$ since they are i.i.d. by assumption. 
The values $\mathbb{P}(\gamma_{i}=1) $ that satisfy this linear matrix inequality belong in some closed convex set. 
%, and we know by assumption that the case $\mathbb{P}(\gamma_{i,k}=1) =1$ belongs in this set by the assumption that \textbf{$V_i(x_i)$ is a Lyapunov function for the closed loop mode of the system.}. 
As a result there is a minimum value $c_i$, given by the semidefinite program \eqref{eq:compute_c}, such that condition \eqref{eq:LMI} is equivalent to $\mathbb{P}(\gamma_{i}=1) \geq c_i $. \qedsymbol
%%
%\begin{eqnarray}\label{eq:auxiliary_SDP}
	%c_i = \sup_{y \in \reals^{n_i}, y \neq 0} \frac{ y^T ( A_{o,i}^T P_i A_{o,i} - \rho_i P_i ) y }
	%{ y^T ( A_{o,i}^T P_i A_{o,i} - A_{c,i}^T P_i A_{c,i} ) y }. 
%\end{eqnarray}
%%
%This is equivalent to the semidefinite program \eqref{eq:compute_c}. \qedsymbol
%%
\end{pf*}

%This proposition transforms control performance of system $i$ to communication variables, i.e., the packet success rate on link $i$. 
The interpretation of the quadratic function $V_i(x_i)$ in this proposition is that it acts as a Lyapunov function for the control system -- see Remark~\ref{rem:Lyapunov}. When the loop closes, the Lyapunov function of the system state decreases, while in open loop it increases, and \eqref{eq:Lyapunov_condition} describes an overall decrease in expectation over the packet success. 
We also point out that an implicit assumption for Proposition~\ref{prop:Lyapunov_abstraction} and throughout the paper is that if system $i$ always remains in closed loop then the desired decrease rate is met, that is, $A_{c,i}^T P_i A_{c,i} \preceq \rho_i P_i$. %The required packet success rate computed in \eqref{eq:compute_c} satisfies $c_i \leq 1$ as long as $A_{c,i}^T P_i A_{c,i} \preceq \rho_i P_i$, i.e., $V_i(x_i)$ is chosen as a Lyapunov function for the stable mode $A_{c,i}$ in \eqref{eq:system}. %Further details about the choice of this abstraction can be found in Remark~\ref{rem:Lyapunov}. 

In this paper we assume that quadratic Lyapunov functions $V_i(x_i)$ and desired expected decrease rates $\rho_i$ are given for each control system. They present a control interface for communication design over a shared wireless medium. We design the sensor access policies so that \textit{all Lyapunov functions $i$} decrease at the desired rates $\rho_i <1$ at any time $k$ \textit{in expectation}. By the above proposition, these control performance requirements are equivalent to minimum packet success rates \eqref{eq:gamma_requirement} for each link $i$, computed by \eqref{eq:compute_c}. Hence we need to ensure that \eqref{eq:gamma_requirement} holds for all links $i$. 

Besides control performance, it is desired that the sensors' channel access mechanism makes an efficient use of their power resources. We assume that when sensor $i$ decides to access the channel it transmits with a fixed power $p_i>0$. The total expected power consumption at each slot is given by $\sum_{i=1}^m \mathbb{E}_{h_i} \alpha_{i}(h_i) p_{i} $ summing up the transmit power of each system $i$ if the system decides to transmit. We pose then the design of the sensor access policies $\alpha$ that minimize the total expected power consumption subject to the desired control performance \eqref{eq:Lyapunov_condition} (equivalently \eqref{eq:gamma_requirement}) for all plants as
\begin{eqnarray}\label{eq:main_optimization}
	\underset{ \alpha \in  \mathcal{A}}  %[0,1]^m
	{\text{minimize}\,\,} & & \sum_{i=1}^m \mathbb{E}_{h_i} \alpha_{i}(h_i) p_{i}  \\
	\text{subject to} & &
	%\mathbb{E}_h \alpha_i(h) \, q(h_i, p_i(h))  \, R_i \succeq R_i' , \; i=1,\ldots,m \notag\\
	c_i \leq \mathbb{E}_{h_i} [ \alpha_i(h_i) q(h_i)]\, \prod_{j\neq i} \Big[ 1 - \mathbb{E}_{h_j} [\alpha_j(h_j)] q_{ji} \Big] \notag\\
	& & \qquad i = 1,\ldots,m. %1 \leq i \leq m.  
	\label{eq:main_constraint}%\\
	%& & &\alpha: \mathcal{H}^m \mapsto \Delta^m, \, p: \mathcal{H}^m \mapsto [0, p_{\max}]^m \notag
	%
\end{eqnarray}
Technically we assume that the problem is strictly feasible. This is a common constraint qualification assumption that will allow us to examine the Lagrange dual problem~\cite[Ch. 5]{Boyd_convex}. More specifically we assume the following.

\begin{asmptn}\label{as:strictly_feasible}
There exists $\alpha' \in \mathcal{A}$ that satisfies constraints \eqref{eq:main_constraint} with strict inequality.
%, i.e., for all $i=1,\ldots,m$
%%
%\begin{eqnarray}\label{eq:strictly_feasible}
	%c_i < \mathbb{E}_{h_i} [ \alpha'_i(h_i) q(h_i)]\, \prod_{j\neq i} \Big[ 1 - \mathbb{E}_{h_j} [\alpha'_j(h_j)] q_{ji} \Big].
%\end{eqnarray}
%%
\end{asmptn}

In the following section we proceed to characterize the optimal access policies $\alpha^*$. 
%, by transforming the original problem \eqref{eq:main_optimization}-\eqref{eq:main_constraint} into an equivalent one that has a zero duality gap. 
In particular we reveal a simple and intuitive decoupled structure. Each sensor $i$ independently accesses the channel in a way that trades off the goal of closed loop $i$ with the effect of collisions on all other closed loops $j \neq i$ collectively. Later in Section~\ref{sec:optimization} we develop a procedure to find these optimal access policies.

\begin{remark} \label{rem:Lyapunov}
In this paper we are interested in communication design for control performance, in contrast to determining what communication designs guarantee stability commonly examined in the literature, e.g.~\cite{Branicky_stability, Hristu_shared_feedback, NCS_Aloha, Hespanha_survey}. The Lyapunov-like abstraction \eqref{eq:Lyapunov_condition} of Proposition~\ref{prop:Lyapunov_abstraction} provides a characterization of control performance, which also implies stability. If \eqref{eq:Lyapunov_condition} holds for each time step $k=0,\ldots,N-1$, then by taking the expectation at both sides and by iterating backwards in time we find that 
\begin{eqnarray}\label{eq:performance_over_time}
	&\mathbb{E} V_i( x_{i,N}) \leq \rho_i
	\, \mathbb{E} V_i(x_{i,N-1}) +  Tr(P_i W_i)  \notag\\
	&\leq \ldots \leq
	\rho_i^N \, \mathbb{E} V_i(x_{i,0}) + %\frac{1-\rho_i^N}{1-\rho_i} \, Tr(P_i W_i)
	\sum_{k=0}^{N-1} \rho_i^k \; Tr(P_i W_i) .
\end{eqnarray}
Hence, system states have second moments that decay exponentially with rate $\rho_i <1$ and in the limit remain bounded by $Tr(P_i W_i) /(1-\rho_i)$, since the sum in \eqref{eq:performance_over_time} converges. As a technical sidenote, an advantage of the Lyapunov performance approach is that it defines a convex region (a lower bound) for the packet success rate in \eqref{eq:gamma_requirement}, which is easy to employ in our random access design \eqref{eq:main_constraint}. %, and thus presents a suitable interface between control and communication design. 
On the contrary, a jump linear system of the form~\eqref{eq:system} is (mean square) stable if and only if the spectral radius of the matrix $\mathbb{P}(\gamma_{i}=1) A_{c,i} \otimes A_{c,i} + \mathbb{P}(\gamma_{i}=0) A_{o,i} \otimes A_{o,i}$ is less than $1$~\cite{MJLS_stability}. However the spectral radius of a non-symmetric matrix is not convex in general, hence it is unclear how to best examine stability in our random access framework. 
\qedopensymbol
\end{remark}

\section{Channel-aware Random Access Design}\label{sec:solution}

Our main result is the following characterization of the optimal access policies for the sensors.

\begin{thm}[Optimal sensor access]\label{thm:optimal_characterization}
Consider a random access architecture with $m$ control loops of the form \eqref{eq:system}, communication modeled by \ref{eq:gamma_probability}, and control performance abstracted by \eqref{eq:Lyapunov_condition}-\eqref{eq:gamma_requirement} for each loop $i=1,\ldots,m$. Consider the design of optimal sensor access policies \eqref{eq:main_optimization}-\eqref{eq:main_constraint}, and let Assumptions \ref{as:non_atomic}, \ref{as:strictly_feasible} hold. Then there exists a matrix of non-negative elements $\nu^* \in \reals_+^{m \times m}$ such that the optimal sensor access policy for each sensor $i = 1,\ldots, m$ is written as
\begin{eqnarray}\label{eq:optimal_access}
	\alpha_i^*(h_i) = \left\{ \begin{array}{ll} 
	1 & \text{if } \; \nu^*_{ii} \, q(h_i) \geq \, p_{i}  + \sum_{j\neq i} \nu^*_{ji}  \, q_{ij} \\
	%\frac{\nu^*_{ii}}{p_{i}  + \sum_{j\neq i} \nu^*_{ji}} \,  q(h_i) \geq 1 \\
	0 & \text{otherwise.}
	\end{array} \right.
\end{eqnarray}
% 
%%
%\begin{eqnarray}\label{eq:optimal_access}
	%\alpha_i^*(h_i) = \ind{ \frac{\nu^*_{ii}}{p_{i}  + \sum_{j\neq i} \nu^*_{ji}} \,  q(h_i) \geq 1 }.
%\end{eqnarray}
%% 
%
\end{thm}

We observe the following interesting facts. First note that the optimal policies are deterministic, that is, given current channel conditions each sensor either transmits or not. Second we note that by the assumed strict monotonicity of the packet success function $q(.)$, the optimal sensor access policies in \eqref{eq:optimal_access} are threshold policies. That is, a sensor transmits only when its corresponding channel quality is above some threshold. The intuitive interpretation is that a sensor should attempt to close its loop only when its channel is sufficiently favorable, i.e., the sensor experiences a sufficiently high current packet success rate. 

Third, and more importantly, the optimal policies are decoupled among the sensors. That is because the policy $\alpha_i^*$ (or equivalently the threshold for sensor $i$) in \eqref{eq:optimal_access} only depends on parameters pertinent to system $i$, i.e., its transmit power $p_i$, and the values $\nu^*_{ii}$ and $\sum_{j\neq i} \nu^*_{ji} \, q_{ij}$ which belong in the $i$th column of matrix $\nu^*$. Hence, as long as the matrix $\nu^*$ is available, each sensor can select its optimal channel access policy independently of what the other sensors are trying to achieve. We note that decentralized threshold-based policies have also been shown to be optimal for general wireless communication networks~\cite{Adireddy_Tong, Qin_Berry, Hu_Ribeiro_RA}. The context differs however, since in these works the objective is thoughput-based utility functions in contrast to the packet success rates used for control systems here. 

As we explain in the proof, the matrix $\nu^*$ technically corresponds to the optimal Lagrange multiplier of an appropriately defined problem (cf. \eqref{eq:primal_optimization}-\eqref{eq:third_ineq}). An intuitive alternative interpretation is as follows. We can think of each diagonal term $\nu^*_{ii}$ as the importance of control performance of system $i$, and of each off-diagonal term $\nu^*_{ji}$ as the collision effect that sensor $i$ has on another system $j$. The optimal access policy for sensor $i$ in \eqref{eq:optimal_access}, or equivalently the optimal channel threshold, trades off the requirement on loop $i$ and the collective negative effect ($\sum_{j \neq i} \nu^*_{ji} q_{ij}$) on all other control loops $j\neq i$. That is because a larger value $\nu^*_{ii}$ corresponds to a lower threshold (sensor transmits more often), while a larger value $\sum_{j \neq i} \nu^*_{ji} q_{ij}$ corresponds to a higher threshold (sensor transmits less often). Note also that the latter summands are normalized by the parameters $q_{ij}$, i.e., the probability that sensor $i$ collides with link $i$ when both sensors transmit. Morever, a high transmit power $p_i$ in \eqref{eq:optimal_access} also implies that sensor $i$ should access the channel less often to limit expenditures.

%%%%%%%
%\input{Channel_aware_proof}
%%%%%%%

The decoupled structure of the optimal sensor access policies in Theorem~\ref{thm:optimal_characterization} relies on knowing the values $\nu^*$. In the following section we develop a distributed iterative procedure to obtain the desired $\nu^*$.
%, which is easily implementable in the architecture of \figref{random_access}. In particular the access point/controller maintains variables $\nu$ and communicates them to the sensors via the reverse channel, and the sensors adapt their access rates according to the decoupled form of Theorem~\ref{thm:optimal_characterization}. 

\begin{remark}
In our previous work in \cite{c_GatsisEtal15} we consider simpler random access policies for the sensors, not taking into account channel state information. In particular we consider that at every time $k$ each sensor $i$ randomly and independently transmits with some constant probability $\tilde{\alpha}_i \in [0,1]$ to be designed. Similarly to \eqref{eq:gamma_probability} the probability of successfully closing each loop is given by
%
%\begin{eqnarray}\label{eq:gamma_probability_no_channel}
	$\mathbb{P}( \gamma_{i,k} = 1) = \tilde{\alpha}_i \, q_{ii} \; \prod_{j\neq i} \Big[ 1 - \tilde{\alpha}_j q_{ji} \Big]$.
%\end{eqnarray}
%
It turns out~\cite[Theorem 2]{c_GatsisEtal15} that the optimal access rates $\tilde{\alpha}^*$, i.e., the solution to a problem equivalent to the channel-aware setup in~\eqref{eq:main_optimization}-\eqref{eq:main_constraint} can be expressed as
\begin{eqnarray}\label{eq:optimal_access_no_channel}
	\tilde{\alpha}_i = \frac{\tilde{\nu}_{ii}} { p_i + \sum_{j \neq i} \tilde{\nu}_{ji} \, q_{ij}} 
\end{eqnarray}
for each $i \in  \{1,\ldots,m\}$ for some non-negative matrix $\tilde{\nu} \in \reals_+^{m \times m}$. The matrix $\tilde{\nu}$ here has the same interpretation as the matrix $\nu^*$ of Theorem~\ref{thm:optimal_characterization} but the two matrices are different as they correspond to different problems. Hence we see that for the non-channel-aware case the sensors need to randomize $(0 < \tilde{\alpha}_i <1)$. In contrast, conditioned on channel state information being available the optimal policies for the sensors are deterministic, exploiting favorable channel conditions to transmit. 
\end{remark}

\section{Computation of Channel-aware Random Access Policies}\label{sec:optimization}

In the previous section the optimal sensor access policies that guarantee control performance of all closed loop systems and minimize power expenditures are characterized in terms of some appropriate matrix of values $\nu^*$ (Theorem~\ref{thm:optimal_characterization}). In this section we capitalize on this result and develop an iterative procedure to determine the optimal sensor access policies by computing these values $\nu^*$. The procedure is distributed and easily implementable in the architecture of \figref{random_access}. In particular the common access point/controller is responsible for finding $\nu^*$ and communicates them to the sensors via the reverse channel, so that the sensors do not need to directly coordinate or communicate among themselves.

Technically as we have argued in the proof of Theorem~\ref{thm:optimal_characterization} the values $\nu^*$ are the optimal Lagrange dual variables of an appropriately defined problem (cf.\eqref{eq:primal_optimization}-\eqref{eq:third_ineq}). The iterative procedure presented in Algorithm~\ref{alg:random_access} corresponds mathematically to a dual subgradient algorithm~\cite[Ch. 8]{Bertsekas_convex} to find the optimal dual variables $\nu^*$. Alternatively we can interpret the procedure as a distributed implementation in the wireless control architecture of \figref{random_access} as follows.

At each period $t$ the access point/controller of \figref{random_access} maintains a tentative matrix of values $\nu(t)$. At the beginning of each period, the access point (AP) sends to each sensor $i$ the values $\nu_{ii}(t)$ and $\sum_{j \neq i} \nu_{ji}(t) q_{ij}$ via the reverse channel (Step~\ref{step:reverse_comm}). For the rest of the period $t$ each sensor uses a random access policy $\alpha(h_i,; t)$ as if the received values $\nu(t)$ corresponded to the optimal $\nu^*$ (Step~\ref{step:primal}). Here $\alpha(h_i; t)$ denotes the valuation of the policy during period $t$ at any channel state $h_i \in \reals_+$. Then the AP measures the gap between desired and current control performance of each system during this period and updates the values $\nu(t)$ to $\nu(t+1)$ to prepare for the next period (Step~\ref{step:dual}). To perform this update the AP needs to compute\footnote{Here we assume that even when collisions arise the AP can identify which sensor transmits at each time slot. Hence it can measure the average rate $\mathbb{E}_{h_i} [ \alpha_i(h_i; t)]$ at which each sensor $i$ accesses the channel, as well as the term $\mathbb{E}_{h_i} [ \alpha_i(h_i;  t) \, q(h_i)]$ which is the packet success ratio when only sensor $i$ transmits.
% which by \ref{eq:gamma_probability} equals the ratio $\mathbb{P}( \gamma_{i} = 1) / \prod_{j\neq i} \Big[ 1 - \mathbb{E}_{h_j} [\alpha_j(h_j;t)] \, q_{ji} \Big]$ of the average packet success rate over the ??
} the average transmission and packet success rates for each system during this period (Step~\ref{step:measure}) and keep track of some auxiliary variables (Step~\ref{step:beta}). 

\begin{algorithm}[!t]{\small
\caption{Distributed random access computation}\label{alg:random_access}
\begin{algorithmic}[1]
%
%\Require ...
%\Statex
\State Initialize $\lambda(0) \in \reals_+^m$, $\nu(0) \in \reals_+^{m \times m}$ at the AP %, $t \gets 0$
\Loop $\quad$ At period $t = 0, 1, \ldots$
\State\label{step:reverse_comm} AP sends $\nu_{ii}(t)$, $\sum_{j \neq i} \nu_{ji}(t) \, q_{ij}$ to each sensor $i$.
\State\label{step:primal} During the period each sensor $i$ accesses the channel according to policy
\begin{eqnarray}\label{eq:primal_step_alpha}
\alpha_i(h_i; t) \gets \left\{ \begin{array}{ll}  1 & \text{if }\nu_{ii}(t) \, q(h_i) \geq p_{i} 
 + \sum_{j\neq i} \nu_{ji}(t) \, q_{ij}\\
0 &\text{otherwise.} 
\end{array}
\right.
\end{eqnarray}
\State\label{step:measure} AP measures $\mathbb{E}_{h_i} [ \alpha_i(h_i;  t) \, q(h_i)]$, $\mathbb{E}_{h_i} [\alpha_i(h_i; t)]$ for all sensors $i=1,\ldots,m$ during the period.
\State\label{step:beta} AP computes the auxiliary variables
\begin{eqnarray}
	\beta_{ii}(t) &\gets &\left[ \frac{\lambda_i(t) }{ \nu_{ii}(t) } \right]_{\mathcal{B}} 
	\label{eq:primal_step_beta_diagonal}\\
	\beta_{ji}(t) &\gets &\left[ 1  - \frac{\lambda_i(t) } { \nu_{ij}(t) } \right]_{\mathcal{B}}
	\label{eq:primal_step_beta_offdiagonal}
\end{eqnarray}
for all $i \neq j \in \{1,\ldots,m\}$, where $[\; ]_{\mathcal{B}}$ denotes the projection to the set defined in \eqref{eq:beta_set_def}.
%
%\State AP computes the matrix $s_\nu(t) \in \reals^{m \times m}$ with diagonal and offdiagonal elements
%%
%\begin{eqnarray}\label{eq:diagonal_subgradient}
	%s_{\nu, ii}(t) &\gets & \beta_{ii}(t) - \mathbb{E}_{h_i} [ \alpha_i(h_i; t) \, q(h_i)] \\
	%%
	%\label{eq:offdiagonal_subgradient}
	%s_{\nu, ij}(t) &\gets & \mathbb{E}_{h_j} [\alpha_j(h_j; t)] - \beta_{ji}(t)
%\end{eqnarray}
%%
%respectively for all $i \neq j \in \{1, \ldots, m\}$, and also the vector $s_\lambda(t) \in \reals^{m }$ with elements
%%
%\begin{eqnarray}\label{eq:vector_subgradient}
	%s_{\lambda, i}(t) \gets \log(c_i) &- & \log(  \beta_{ii}(t)) \, \notag\\
	%&- & \sum_{j\neq i} \log ( 1 - \beta_{ji}(t) q_{ji} ) 
%\end{eqnarray}
%%
%%
%\State AP updates the new dual variables
%%
%\begin{eqnarray}
	%\nu (t+1) \gets \Big[ \nu(t) + \epsilon(t) \, s_\nu (t) \Big]_+   \label{eq:dual_step_nu} \\
	%\lambda (t+1) \gets \Big[ \lambda(t) + \epsilon(t) \, s_\lambda(t) \Big]_+   \label{eq:dual_step_lambda}
%\end{eqnarray}
%%
%where $[\; ]_+$ denotes the elementwise projection to the non-negatives $\reals_+^{m \times m}$.
%%
\State\label{step:dual} AP updates the new dual variables
\begin{eqnarray}
	&\nu_{ii} (t+1) \gets &\Big[ \nu_{ii}(t) +  \notag\\
	&&
	\epsilon(t) \,
	\big( \beta_{ii}(t) - \mathbb{E}_{h_i} [ \alpha_i(h_i; t)  q(h_i)] \big)  \Big]_+   \label{eq:dual_step_nu_diag} \\
	&\nu_{ij} (t+1) \gets &\Big[ \nu_{ij}(t) + \notag\\
	&&\epsilon(t) \, 
	\big( \, \mathbb{E}_{h_j} [\alpha_j(h_j; t)] q_{ji} - \beta_{ji}(t) \, \big) \, \Big]_+  	\label{eq:dual_step_nu_offdiag}\\
	&\lambda_i (t+1) \gets &\Big[ \lambda_i(t) + \epsilon(t) \, \big( \, \log(c_i) -  \log(  \beta_{ii}(t)) \, \notag\\
	&&-  \sum_{j\neq i} \log ( 1 - \beta_{ji}(t) ) \,  \big) \, \Big]_+   \label{eq:dual_step_lambda}
\end{eqnarray}
for all $i \neq j \in \{1, \ldots, m\}$, where $[\; ]_+$ denotes the projection to the non-negatives $\reals_+$.
%
%\State $k \gets k+1 $
%
\EndLoop
\end{algorithmic}}
\end{algorithm}

This algorithm is guaranteed to converge to the optimal sensor access policies as we state next.

\begin{thm}[Sensor access optimization]\label{thm:convergence}
Consider the setup of Theorem~\ref{thm:optimal_characterization}. 
%Consoder the design of optimal sensor access policies in \eqref{eq:main_optimization}-\eqref{eq:main_constraint}, and let Assumptions \ref{as:non_atomic}, \ref{as:strictly_feasible} hold. 
The iterations of Algorithm~\ref{alg:random_access} with stepsizes in \eqref{eq:dual_step_nu_diag}-\eqref{eq:dual_step_lambda} satisfying $\sum_{t\geq 0} \epsilon(t)^2 < \infty$, $\sum_{t \geq 0} \epsilon(t) = \infty$ converge to the optimal sensor access policies, i.e.,
\begin{eqnarray} \label{eq:alg_constraint_conv}
	&&c_i \leq \lim_{t \rightarrow \infty} 
	\mathbb{E}_{h_i} [ \alpha_i(h_i;t) q(h_i)]\, \prod_{j\neq i} \Big[ 1 - \mathbb{E}_{h_j} [\alpha_j(h_j;t)] q_{ji} \Big] ,
	\notag \\
	&&
\end{eqnarray}
for all $i=1,\ldots,m$, and 
\begin{eqnarray} \label{eq:alg_value_conv}
	\lim_{t \rightarrow \infty} \, \sum_{i=1}^m \mathbb{E}_{h_i} [\alpha_{i}(h_i;t)] p_{i} 
	= \sum_{i=1}^m \mathbb{E}_{h_i} [\alpha_{i}^*(h_i)] p_{i} .
\end{eqnarray}
\end{thm}

%Apart from converging to the optimal operating point, Algorithm~\ref{alg:random_access} is easily implementable in the wireless control architecture of \figref{random_access}. The sensors decide upon their access rates without coordination among themselves. Moreover they do not need to know global problem information, e.g., the specifications of the other coexisting control loops, or even how many other sensors are sharing the same wireless medium. Each sensor only needs to know the amount of collisions it causes on all other sensors collectively (captured by the value of the sum $\sum_{j \neq i} \nu_{ji}(t)$ in \eqref{eq:primal_step_alpha}).

The caveat of this distributed implementation is that it requires information exchange between sensors and the access point, hence it introduces some communication overhead. This overhead however burdens mainly the access point which is typically a base station with more capabilities compared to the simpler wireless sensors.

%\begin{remark}\label{rem:noninterfering}
%%
%The problem formulation can be modified to include the case where some sensors are not interfering with each other. If sensor $j$ never causes collisions to transmission $i$ we have $q_{ji}=0$. Define the subset of sensors that affect link $i$ as $I_i = \{ j \neq i : q_{ji} >0 \}$, and conversely the subset of links that are affected by sensor $i$ as $O_i = \{ j \neq i: q_{ij} >0\} $. Then the packet success probability in \eqref{eq:gamma_probability} is modified so that the product is over the interfering sensors $\prod_{j \in I_i}$. Similarly the optimal sensor access rates in \eqref{eq:optimal_access} are modified to include the sum $\sum_{j \in O_i } \nu^*_{ji}$. Algorithm \ref{alg:random_access} is also modified so that no 'coupling' variables $\alpha_{ji}, \nu_{ij}$ are needed when $j \notin I_i$.
%%
%\qedopensymbol
%\end{remark}

\section{Numerical simulations}\label{sec:simulations}

\begin{figure}[!t]
\centering
\includegraphics[width=\columnwidth]{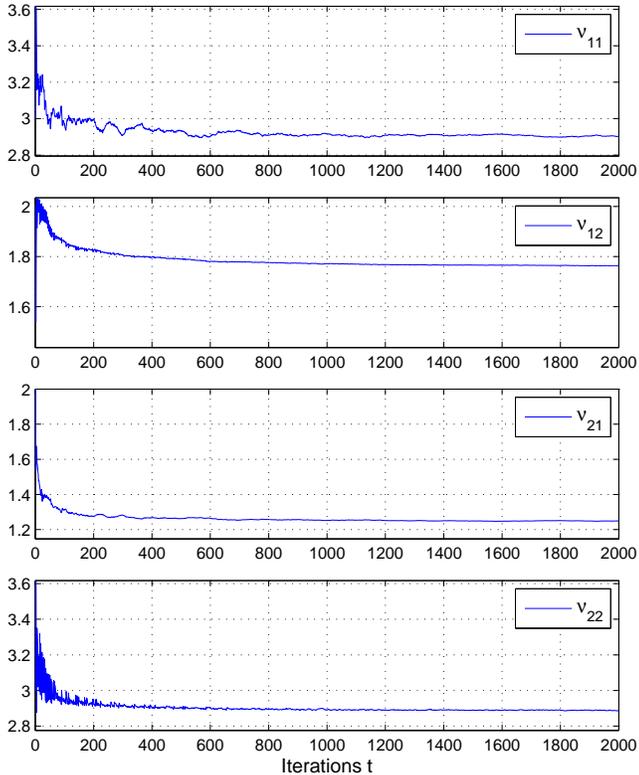}
\caption{Evolution of dual variables during the optimization algorithm. The elements of the matrix $\nu(t)$ converge to the optimal values $\nu^*$ required to obtain the optimal sensor access policies.}
\label{fig:dual_evolution}
\end{figure}

\begin{figure}[!t]
\centering
\includegraphics[width=\columnwidth]{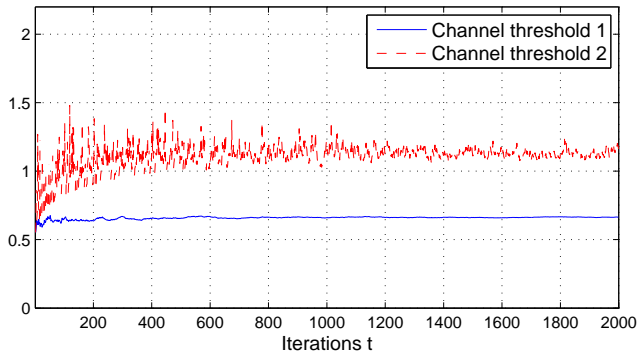}
\caption{Channel thresholds corresponding to the access policies selected by the optimization algorithm. The channel thresholds for both sensors converge to their optimal values in the limit. Sensor 1 has a lower threshold, i.e., transmits more often, since it is required to guarantee control performance for a more demanding (unstable) plant.}
\label{fig:channel_thresholds}
\end{figure}

\begin{figure}[!t]
\centering
\includegraphics[width=\columnwidth]{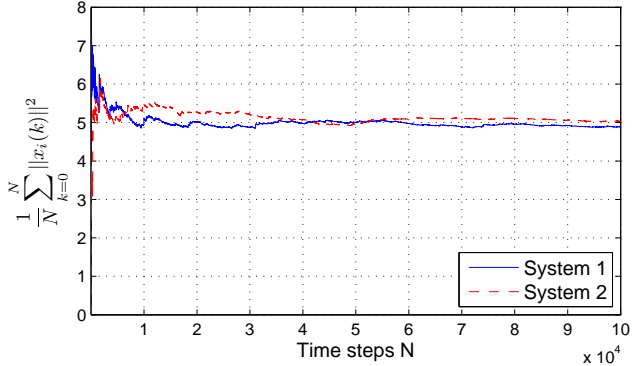}
\caption{Evolution of control systems using the optimal random access policies. Both systems remain stable despite collisions and packet drops. Also their long run average quadratic cost converges to the same value, since be design both systems were required to have the same control performance.}
\label{fig:systems_simulation}
\end{figure}

We present a numerical example of the random access design. We consider a case with $m=2$ scalar control systems of the form \eqref{eq:system}. We assume the first system has open and closed loop dynamics given by $A_{o,1} = 1.1$, $A_{c,1} = 0.5$ respectively, i.e., it is open loop unstable. We assume the second system has integrator open loop dynamics $A_{o,i} = 1$ and stable closed loop dynamics $A_{c,2} = 0.4$. Both systems are perturbed by zero-mean unit-variance Gaussian noises, hence both system states will diverge unless the closed loops are applied appropriately. The systems are asymmetric, but we model a symmetric control performance requirement. The Lyapunov function $V_i(x_i) = x_i^2$ ($P_i =1$) for both plants $i=1,2$ is required to decrease with expected rate $\rho_1 = \rho_2 = 0.8$ (cf.~\eqref{eq:Lyapunov_condition}). By Proposition~\ref{prop:Lyapunov_abstraction} these control performance requirements are equivalent to required packet success rates $c_1 \approx 0.43$, $c_2 \approx 0.27$ for the two sensors, computed by \eqref{eq:compute_c}. Hence System 1, which is more unstable, requires a higher packet success rate.

We assume that both channel states $h_{1,k}, h_{2,k}$ are i.i.d. exponential with mean $1$. In isolation each sensor faces a packet success probability modeled by the function $q(h_{i,k})$, $i=1,2$ shown in \figref{q_vs_SNR}. Also when both sensors transmit at the same time, collisions occur with probability $q_{12}=q_{21} = 0.5$. The transmit powers are taken equal $p_1 = p_2 = 1$. 

We solve the random access design problem~\eqref{eq:main_optimization}-\eqref{eq:main_constraint} by implementing Algorithm~\ref{alg:random_access}, which as explained in Section~\ref{sec:optimization} solves the problem in the dual domain. We note that at each iteration of the algorithm some expectations with respect to the channel state distributions need to be computed, in particular in steps \ref{eq:dual_step_nu_diag} and \ref{eq:dual_step_nu_offdiag} of Algorithm~\ref{alg:random_access}. In our simulations we approximate these expectations with averages from a large number of samples, since samples from the exponential channel distributions can be readily simulated. The iterates of the matrix dual variables $\nu(t)$ during the simulation are shown in \figref{dual_evolution} where we observe that they converge to the optimal values $\nu^*$, as was also shown in the proof of Theorem~\ref{thm:convergence}. We also plot the evolution of the sensor access policies $\alpha_i(h_i; t)$, or equivalently the thresholds of these policies during the simulation of the algorithm in \figref{channel_thresholds}. As also established in Theorem~\ref{thm:convergence} the channel thresholds converge to their optimal values in the limit. We observe that the Sensor 1 has a lower threshold, meaning that it transmits more often, which is natural since it corresponds to the unstable plant. 

After the optimal access policies (equivalently channel thresholds) have been found, we simulate the random access architecture with the obtained. In \figref{systems_simulation} we plot the empirical average long term quadratic cost of the systems $1/N \sum_{k=1}^N x_{i,k}^2$ for each system $i=1,2$. We first observe that both systems remain stable despite packet collisions over the shared channel. Moreover, even though the two systems are asymmetric, both long term average costs converge to the same value because we required the same control performance for both systems. More specifically this long term cost equals the value $\text{Tr}(P_i W_i )/(1 - \rho_i ) = 1/(1-0.8) = 5$ for both systems $i=1,2$, as noted in Remark~\ref{rem:Lyapunov}. Hence even though the two plants have different dynamics, the obtained channel access policies provide symmetric performance by design. The empirical rates $1/N \sum_{k=1}^N \alpha_{i,k}$ at which each sensor transmits equal $0.51$ and $0.32$ for $i=1,2$ respectively. As expected, both sensors access the channel at a rate higher than the respective necessary packet success rate on each link, i.e., $\alpha_i^* > c_i$. This happens because the sensors need to counteract the effect of packet collisions, as well as packet drops due to decoding errors. 

% Other things I could include: 
%- the rates at which a sensor transmits $\mathbb{E} \alpha$, which is larger for system 1. 
%- Perhaps the term $\mathbb{E} \alpha q$ but this is not so intuitive 
%- The cost/power in the limit as it is minimized?

%Comparisons with my channel-aware centralized scheduling??
%%
\section{Concluding Remarks}\label{sec:conclusion}

We design a random access mechanism for sensors transmitting measurements of multiple plants over a shared wireless channel to a controller. The goal of the sensors is to guarantee control performance for all control systems by mitigating the effect of packet collisions from simultaneous transmissions as well as by adapting online to randomly varying channel conditions. Via a Lyapunov function abstraction, control performance is transformed to required packet success rates of each closed loop. We show that the optimal random access policies can be decoupled between the sensors and are of a threshold form with respect to channel states. Moreover we develop a distributed procedure to compute the optimal policies.

In future research we aim to explore the ability of the distributed random access design procedure to track changes in the problem parameters, for example, variations in the channel collision pattern, changes in the control performance requirements, or the admission of new control loops in the architecture. Additionally some preliminary work on online sensor adaptation to plant states as in, e.g., the single link case of~\cite{GatsisEtal14} or the scheduling of~\cite{Walsh_stability, Donkers_switched}, is considered in \cite{ICCPS16_GatsisEtal}.

\bibliography{control_multiple_access}           % and a bib file to produce the 
                                 % bibliography (preferred). The
                                 % correct style is generated by
                                 % Elsevier at the time of printing.

\appendix

\section{Proof of Theorem~\ref{thm:optimal_characterization}}

%\begin{pf*}{Proof of Theorem~\ref{thm:optimal_characterization}.}
%
The first part of the proof involves converting problem \eqref{eq:main_optimization}-\eqref{eq:main_constraint} into an equivalent auxiliary optimization problem which has zero duality gap. Then in the second part we use Lagrange duality arguments to show that \eqref{eq:optimal_access} describes an optimal solution for the auxiliary problem.

We begin by a modification to remove the product of the expectations appearing in the constraints \eqref{eq:main_constraint}. Taking the logarithm at each side of \eqref{eq:main_constraint} preserves the feasible set of variables by monotonicity. Then the logarithm of the product at the right hand side of \eqref{eq:main_constraint} becomes a sum of logarithms, and we can rewrite the optimal random access design problem equivalently as %\KG{maybe use the normal align command and force from the beginning all aligns left??}
\begin{eqnarray}\label{eq:main_optimization_1}
	\underset{ \alpha \in \mathcal{A} }
	{\text{minimize}\,\,} &\; & \sum_{i=1}^m \mathbb{E}_{h_i} \alpha_{i}(h_i) p_{i}  \\
	\text{subject to}  
	%\mathbb{E}_h \alpha_i(h) \, q(h_i, p_i(h))  \, R_i \succeq R_i' , \; i=1,\ldots,m \notag\\
	&\; &\log(c_i) \leq \log( \mathbb{E}_{h_i} [ \alpha_i(h_i) \, q(h_i)] ) \notag\\
	& & \qquad + \, \sum_{j\neq i} \log ( 1 - \mathbb{E}_{h_j} [\alpha_j(h_j)]  q_{ji} ),
	\label{eq:main_constraint_1} \\
	& &i=1,\ldots,m.  \notag%\\
	%& & &\alpha: \mathcal{H}^m \mapsto \Delta^m, \, p: \mathcal{H}^m \mapsto [0, p_{\max}]^m \notag
	%
\end{eqnarray}
Here we make an implicit technical assumption that the terms $\mathbb{E}_{h_i} [ \alpha_i(h_i) \, q(h_i)]$ and $\mathbb{E}_{h_i} [ \alpha_i(h_i) ] \, q_{ij}$ in \eqref{eq:main_constraint_1}, which in general take values in the unit interval $[0,1]$ as all involved variables belong there too, are bounded away from $0$ and $1$. Then the logarithms in \eqref{eq:main_constraint_1} are well-defined and finite. This does not restrict the feasible set of solutions, as intuitively each sensor $i$ can neither choose $\alpha_i(h_i)$ too close to 0 otherwise it cannot meet its packet success requirement in \eqref{eq:main_constraint}, nor too close to 1 otherwise it causes significant packet collisions on other sensors. 

%Next, we introduce auxiliary variables in place of the expectations in the constraint \eqref{eq:main_constraint_1} for $i=1,\ldots,m$. In particular we replace the term $\mathbb{E}_{h_i} [ \alpha_i(h_i) \, q(h_i)]$ by a variable $\beta_{ii}$, and the terms $\mathbb{E}_{h_j} [\alpha_j(h_j)]$ by variables $\beta_{ji}$ for $j \neq i$. Hence we rewrite \eqref{eq:main_constraint_1} as

Next, we replace the term $\mathbb{E}_{h_i} [ \alpha_i(h_i) \, q(h_i)]$ in constraint \eqref{eq:main_constraint_1} by an auxiliary variable $\beta_{ii}$ for $i=1,\ldots,m$, and the terms $\mathbb{E}_{h_j} [\alpha_j(h_j)] q_{ji}$ in \eqref{eq:main_constraint_1} by variables $\beta_{ji}$ for $j \neq i$. Hence we rewrite \eqref{eq:main_constraint_1} as
\begin{eqnarray}\label{eq:aux_constraint}
	\log(c_i) \leq \log(  \beta_{ii}) \, + \sum_{j\neq i} \log ( 1 - \beta_{ji} ).
\end{eqnarray}
To force the auxiliary variables to behave like the expectations we introduce additional constraints of the form  
\begin{eqnarray}\label{eq:aux_constraint_2}
	&\beta_{ii} &\leq \mathbb{E}_{h_i} [ \alpha_i(h_i) \, q(h_i)] \\
	&\beta_{ji} &\geq  \mathbb{E}_{h_j} [\alpha_j(h_j)] \, q_{ji}
\end{eqnarray}
for all $i, j \in \{1,\ldots,m \} , \, j \neq i$. Each of these variables are restricted to a subset
\begin{eqnarray}\label{eq:beta_set_def}
	\beta_{ij} \in \mathcal{B} = [\beta_{\min} , \beta_{\max} ]
\end{eqnarray}
of the unit interval $[0,1]$. In a matrix form $\beta \in \mathcal{B}^{m \times m}$. These upper and lower bounds guarantee that all logarithms at constraints \eqref{eq:aux_constraint} are finite, as we also assumed for constraint \eqref{eq:main_constraint_1}. Overall we formulate the auxiliary optimization problem
%%
%\begin{align}
	%&\underset{ \alpha \in \mathcal{A}, \, \beta \in \mathcal{B} }
	%{\text{minimize}} & & \sum_{i=1}^m \mathbb{E}_{h_i} \alpha_{i}(h_i) p_{i} \label{eq:primal_optimization}\\
	%%
	%&\text{subject to} & &
	%%\mathbb{E}_h \alpha_i(h) \, q(h_i, p_i(h))  \, R_i \succeq R_i' , \; i=1,\ldots,m \notag\\
	%\log(c_i) \leq \log(  \beta_{ii}) \, + \sum_{j\neq i} \log ( 1 - \beta_{ji} q_{ji} )  \label{eq:first_ineq} \\
	%& & & \beta_{ii} \leq \mathbb{E}_{h_i} [ \alpha_i(h_i) \, q(h_i)] \label{eq:second_ineq} \\
	%& & & \beta_{ji} \geq  \mathbb{E}_{h_j} [\alpha_j(h_j)]  \label{eq:third_ineq} \\
	%& & &i, j \in \{1,\ldots,m \} , \, j \neq i \notag
	%%
%\end{align}
%%
%%%%%%%%% TRY WITH EQNARRAY
\begin{eqnarray}
	&\underset{ \alpha \in \mathcal{A}, \, \beta \in \mathcal{B}^{m \times m} }
	{\text{minimize}\,\,}  \;& \sum_{i=1}^m \mathbb{E}_{h_i} \alpha_{i}(h_i) p_{i}  \label{eq:primal_optimization}\\
	&\text{subject to} 
	\;& \log(c_i) \leq \log(  \beta_{ii}) \, + \sum_{j\neq i} \log ( 1 - \beta_{ji} )  \label{eq:first_ineq} \\
	& &  \beta_{ii} \leq \mathbb{E}_{h_i} [ \alpha_i(h_i) \, q(h_i)] \label{eq:second_ineq} \\
	& &  \beta_{ji} \geq  \mathbb{E}_{h_j} [\alpha_j(h_j)] \,  q_{ji} \label{eq:third_ineq} \\
	& & i, j \in \{1,\ldots,m \} , \, j \neq i \notag
\end{eqnarray}
We argue that this auxiliary problem is equivalent to the original one in \eqref{eq:main_optimization}-\eqref{eq:main_constraint}, in the sense that a feasible solution of one problem corresponds to a feasible solution with the same objective value for the other problem. Indeed let's start with a feasible solution $\alpha$ for \eqref{eq:main_optimization}-\eqref{eq:main_constraint}. Let us define variables $\beta$ that make \eqref{eq:second_ineq}, \eqref{eq:third_ineq} hold with equality for all $i, j \in \{1,\ldots,m \} , \, j \neq i$. Then the pair $\alpha, \beta$ is also feasible for problem \eqref{eq:primal_optimization}-\eqref{eq:third_ineq} and has the same objective. Reversely, consider a feasible pair $\alpha, \beta$ for problem \eqref{eq:primal_optimization}-\eqref{eq:third_ineq}. Without loss we can assume that all constraints \eqref{eq:second_ineq}-\eqref{eq:third_ineq} hold with equality. Otherwise if, say, an inequality $i$ in \eqref{eq:second_ineq} is strict, we can increase the value of variable $\beta_{ii}$ till equality in \eqref{eq:second_ineq} is reached without loss of feasibility in \eqref{eq:first_ineq} and without changing the objective value in \eqref{eq:primal_optimization}. A similar procedure can be performed if some inequality \eqref{eq:third_ineq} is strict, leading finally to a new feasible point satisfying \eqref{eq:second_ineq}-\eqref{eq:third_ineq} with equalities. Then it is immediate that $\alpha$ is also feasible for \eqref{eq:main_constraint} and has the same objective.

Based on the established equivalence, in the rest of the proof it suffices to show that \eqref{eq:optimal_access} describes an optimal solution for the auxiliary problem \eqref{eq:primal_optimization}-\eqref{eq:third_ineq}. The advantage of formulating this auxiliary problem is that it has zero duality gap as can be shown by the results in \cite{Ribeiro_EURASIP}. To formally state this result, let us denote the optimal value of this problem by $P^*$ (finite by feasibility Assumption~\ref{as:strictly_feasible}) and let us define the Lagrange dual problem. We associate dual variables $\lambda_{i} \geq 0$ with inequalities \eqref{eq:first_ineq}, $\nu_{ii} \geq 0$ with \eqref{eq:second_ineq}, and $\nu_{ij} \geq 0$ with \eqref{eq:third_ineq}, for $i,j \in \{1,\ldots,m\}$. We write the Lagrangian function as
\begin{eqnarray}\label{eq:Lagrangian}
	&&L(\alpha, \beta, \lambda, \nu) = \sum_{i=1}^m \mathbb{E}_{h_i} \alpha_{i}(h_i) p_{i}  
	\notag\\
	&& + \sum_{i=1}^m \lambda_i \Big\{ \log(c_i) - \log(  \beta_{ii}) - \sum_{j\neq i} \log ( 1 - \beta_{ji} ) \Big\}  \notag\\
	&& + \sum_{i=1}^m \nu_{ii} ( \beta_{ii} - \mathbb{E}_{h_i} [ \alpha_i(h_i) \, q(h_i)] )  \notag\\
	&&+ \sum_{i=1}^m \sum_{j\neq i} \nu_{ij} ( \mathbb{E}_{h_j} [\alpha_j(h_j)] \, q_{ji} - \beta_{ji} ) \quad 
\end{eqnarray}
%
%%
%\begin{eqnarray}\label{eq:Lagrangian}
	%&&L(\alpha, \beta, \lambda, \nu) = \sum_{i=1}^m \; \mathbb{E}_{h_i} \alpha_{i}(h_i) p_{i}  \notag\\
	%%
	%&& \qquad+ \sum_{i=1}^m \; \lambda_i ( \log(c_i) - \log(  \beta_{ii}) - \sum_{j\neq i} \log ( 1 - \beta_{ji} q_{ji} ) )  \notag\\
	%&& \qquad + \sum_{i=1}^m \; \nu_{ii} ( \beta_{ii} - \mathbb{E}_{h_i} [ \alpha_i(h_i) \, q(h_i)] ) \notag\\
	%&& \qquad + \sum_{i=1}^m \; \sum_{j\neq i} \nu_{ij} ( \mathbb{E}_{h_j} [\alpha_j(h_j)] - \beta_{ji} ) \quad 
%\end{eqnarray}
%%
Here the dual variables take values $\lambda \in \reals_+^m$, $\nu \in \reals_+^{m \times m}$. We can rearrange the terms of the Lagrangian in the form
\begin{eqnarray}\label{eq:Lagrangian_separable}
	&&L(\alpha, \beta, \lambda, \nu) = \sum_{i=1}^m \Big\{ 
	\mathbb{E}_{h_i} \, \alpha_{i}(h_i)  \, \Big[ p_{i}  + \sum_{j\neq i} \nu_{ji} q_{ij}- \nu_{ii} q(h_i) \Big]  \notag\\
	&& \qquad + \sum_{j\neq i} \left[  - \lambda_i \log ( 1 - \beta_{ji} ) - \nu_{ij} \beta_{ji} \right]\notag\\ 
	&& \qquad + \nu_{ii} \beta_{ii} -   \lambda_i \log(  \beta_{ii})
	+ \lambda_i \log(c_i) \quad \Big\}.
\end{eqnarray}
This form is useful because each primal variable ($\alpha_i(h_i)$ and $\beta_{ji}$ for each $i,j$) is decoupled from the others, a fact we will exploit next. Then we can define the Lagrange dual function
\begin{eqnarray}\label{eq:dual_function}
	g(\lambda, \nu) =  \inf_{\alpha \in \mathcal{A}, \, \beta \in \mathcal{B}^{m \times m} } \; L(\alpha, \beta, \lambda, \nu) ,
\end{eqnarray}
as well as the Lagrange dual problem whose optimal value we denote by $D^*$ as
\begin{eqnarray}\label{eq:dual_problem}
	D^* = \inf_{ \lambda \in \reals_+^m, \, \nu \in \reals_+^{m \times m} } \; g(\lambda, \nu) .
\end{eqnarray}
Then we can establish the following zero duality property about the auxiliary problem \eqref{eq:primal_optimization}-\eqref{eq:third_ineq}.

\begin{proposition}[Strong Duality]
Let Assumptions \ref{as:non_atomic} and \ref{as:strictly_feasible} hold. Then the problem \eqref{eq:primal_optimization}-\eqref{eq:third_ineq} has zero duality gap, i.e., $P^* = D^*$. Moreover if $\alpha^*, \beta^*$ are optimal solutions and $\lambda^*, \nu^*$ are optimal solutions for the dual problem \eqref{eq:dual_problem}, then
\begin{eqnarray}\label{eq:optimal_primal}
	\alpha^*, \beta^* \, \in \argmin_{\alpha \in \mathcal{A}, \, \beta \in \mathcal{B}^{m \times m} } \; L(\alpha, \beta, \lambda^*, \nu^*).
\end{eqnarray}
\end{proposition}

The result follows from \cite[Theorems 1 and 4]{Ribeiro_EURASIP} where general stochastic optimization problems of the form \eqref{eq:primal_optimization}-\eqref{eq:third_ineq} are examined under non-atomic probability measures (Assumption \ref{as:non_atomic}) and strict feasibility (Assumption \ref{as:strictly_feasible}). The proof is omitted due to space limitations. 

The above characterization suggests that we can recover the optimal variables $\alpha^*, \beta^*$ of our problem by just minimizing the unconstrained Lagrangian function. A technical caveat of \eqref{eq:optimal_primal} is that it describes an inclusion only, implying that in general there might be Lagrangian minimizers that are not optimal. The following lemma excludes such cases by establishing that the Lagrangian minimizers $\alpha$, which are functions, i.e., infinite-dimensional variables, are unique up to a set of measure zero. Moreover the following lemma gives an explicit expression for these minimizers.

\begin{lemma}\label{lemma:Lagrange_minimizers}
Consider any dual variables $\lambda \in \reals_+^m, \, \nu \in \reals_+^{m \times m}$. Then the functions $\alpha \in \mathcal{A}$ that minimize the Lagrangian $L(\alpha, \beta, \lambda, \nu)$ are uniquely defined except for a set of arguments $h \in \reals_+^m$ of measure zero, and are given by 
\begin{eqnarray}\label{eq:optimal_alpha}
	\alpha_i(h_i; \lambda, \nu) = \left\{ \begin{array}{ll} 
	1 & \text{if } \; \nu_{ii}q(h_i)\geq p_{i}  + \sum_{j\neq i} \nu_{ji} q_{ij} \\
	0 & \text{otherwise.}
	\end{array} \right.
\end{eqnarray}
for each $i=1, \ldots, m$ and for every value $h_i \in \reals_+$. 
\end{lemma}

In \eqref{eq:optimal_alpha} the term $\alpha_i(h_i; \lambda, \nu)$ denotes the function $\alpha_i$ that minimizes the Lagrangian $L(\alpha, \beta, \lambda, \nu)$ at given dual points $\lambda, \nu$ evaluated at an argument $h_i$. 
%We note that the optimal solution is defined separately for each $i=1,\ldots,m$ and for each argument $h_i$, and also it does not depend on the dual variable $\lambda$ \KG{place elsewhere?}. 
The proof can be found in Appendix~\ref{sec:proof_Lagrangian}.

To sum up we have shown in \eqref{eq:optimal_primal} that the optimal solution $\alpha(.)$ to problem \eqref{eq:primal_optimization}-\eqref{eq:third_ineq} belongs in the set of Lagrange minimizers at $\lambda^*, \nu^*$, and by Lemma~\ref{lemma:Lagrange_minimizers} these minimizers are unique up to a set of measure zero. As a result, all these minimizers will have the same objective and constraint slack in problem \eqref{eq:primal_optimization}-\eqref{eq:third_ineq}, and they will all be optimal for this problem. In particular, the specific minimizer defined by $\alpha_i(h_i; \lambda^*, \nu^*)$ given in \eqref{eq:optimal_alpha} will be optimal for the problem, and corresponds exactly to the one given in \eqref{eq:optimal_access} at the statement of the theorem. 
%
%******************
%EXTRA
%
%\begin{lemma}\label{lemma:Lagrange_minimizers_second?}
%%
%Let Assumption \ref{as:non_atomic} hold and consider any dual variables $\lambda \in \reals_+^m, \, \nu \in \reals_+^{m \times m}$. Then 
%%
%\begin{enumerate}
%%
%\item All pairs $\alpha, \beta \, \in \argmin_{\alpha \in \mathcal{A}, \, \beta \in \mathcal{B} } \; L(\alpha, \beta, \lambda, \nu)$ have the same feasibility slack and objective value with respect to problem \eqref{eq:primal_optimization}-\eqref{eq:third_ineq} i.e., there exist a unique matrix $s_\nu \in \reals^{m \times m}$ and vector vector $s_\lambda \in \reals^{m }$ such that for all minimizers we have
%%
%\begin{eqnarray}%\label{eq:diagonal_subgradient}
	%s_{\nu, ii} &= & \beta_{ii} - \mathbb{E}_{h_i} [ \alpha_i(h_i) \, q(h_i)] \\
	%%
	%%\label{eq:offdiagonal_subgradient}
	%s_{\nu, ij} &= & \mathbb{E}_{h_j} [\alpha_j(h_j)] - \beta_{ji} \\
	%%
	%s_{\lambda, i} &=& \log(c_i) - \log(  \beta_{ii}) - \sum_{j\neq i} \log ( 1 - \beta_{ji} q_{ji} ) 
%\end{eqnarray}
%%
%respectively for all $i \neq j \in \{1, \ldots, m\}$.
%%
%\end{enumerate}
%%
%\end{lemma}
%
%*****************
\qedsymbol

%\end{pf*}

\section{Proof of Lemma~\ref{lemma:Lagrange_minimizers}}\label{sec:proof_Lagrangian}

Consider the problem of minimizing the Lagrangian $L(\alpha, \beta, \lambda, \nu)$ over variables $\alpha \in \mathcal{A}, \beta \in \mathcal{B}^{m \times m}$. Due to the separability of the Lagrangian given in the form \eqref{eq:Lagrangian_separable} over variables $\alpha, \beta$, we can separate the problem into subproblems
\begin{eqnarray}\label{eq:optimal_primal_i}
	\argmin_{\beta_{ji} \in \mathcal{B}} &\; &- \lambda_i \log ( 1 - \beta_{ji} ) - \nu_{ij} \beta_{ji} \\ 
	\argmin_{\beta_{ii} \in \mathcal{B}} &\;  &\nu_{ii} \beta_{ii} -   \lambda_i \log(  \beta_{ii}) 
	\label{eq:optimal_primal_ii}\\
	\argmin_{\alpha_i \in \mathcal{A}_i } &\; &\mathbb{E}_{h_i} \, \alpha_{i}(h_i)  \, \Big[ p_{i}  + \sum_{j\neq i} \nu_{ji} q_{ij} - \nu_{ii} q(h_i) \Big] 
	\label{eq:optimal_primal_iii}
\end{eqnarray}
for $i, j \in \{1, \ldots, m\}, i \neq j$. Next we need to verify that \eqref{eq:optimal_alpha} is optimal for \eqref{eq:optimal_primal_iii}. Note that without loss of generality we can exchange the expectation operator $\mathbb{E}_{h_i}$ and the minimization over $\alpha_i \in \mathcal{A}_i$, which is a function $\alpha_i : \reals_+ \rightarrow [0,1]$ defined for any channel value $h_i \in \reals_+$,  to equivalently solve
\begin{eqnarray}\label{eq:channel_decouple}
	\argmin_{\alpha_i(h_i) \in [0,1] } \; \alpha_{i}(h_i)  \, \Big[ p_{i}  + \sum_{j\neq i} \nu^*_{ji} q_{ij} - \nu^*_{ii} q(h_i) \Big].
\end{eqnarray}
pointwise at all values $h_i \in \reals_+$. 
This is valid because any function $\alpha_i$ that minimizes \eqref{eq:optimal_primal_iii} can differ form the minimizer in \eqref{eq:channel_decouple} at a set of values $h_i \in \reals_+$ with measure at most zero. 

Then we can verify that \eqref{eq:optimal_alpha} is the minimizer in \eqref{eq:channel_decouple}. That is because the right hand side in \eqref{eq:channel_decouple} is a linear expression of $\alpha_i(h_i) \in [0,1]$. Hence the minimizer $\alpha_i(h_i)$ is uniquely defined, and takes values either $0$ or $1$ except for the values $h_i$ where $p_{i}  + \sum_{j\neq i} \nu^*_{ji} q_{ij} - \nu^*_{ii} q(h_i) = 0$. In the latter case the minimizer is not uniquely defined. However due to the strict monotonicity assumption for $q(h_i)$ this case occurs for at most one value $h_i$, hence it is a measure zero event since measure $\phi_i$ is non-atomic by Assumption~\ref{as:non_atomic}. This completes the proof.

%****We note here for future reference a side conclusion. Each one of the summands in the sum in \eqref{eq:Lagrangian_at_optimal} must be zero. This is a complementary slackness property \KG{needed? actually there are a lot of them! for all $i, j \in \{1,\ldots, m\}$, $j \neq i$. }

We also note for future reference the terms $\beta$ that minimize the Lagrangian. Since \eqref{eq:optimal_primal_i}, \eqref{eq:optimal_primal_ii} are strongly convex, their minimizers are unique and satisfy the first order conditions $\partial L/\partial \beta =0$, that is 
\begin{eqnarray}\label{eq:first_order_condition}
	&&\nu_{ii} - \frac{\lambda_i}{\beta_{ii}} 	= 0 \\
	&&\frac{ \lambda_i }{1-\beta_{ji} } - \nu_{ij} 	=0 .
\end{eqnarray}
respectively subject to the box constraints $\beta_{ji} \in \mathcal{B}$ for all $i,j \in \{1,\ldots, m\}$. As a result the optimal solutions are given by
\begin{eqnarray}
	\beta_{ii}(\lambda, \nu) &= &\left[ \frac{\lambda_i }{ \nu_{ii} } \right]_{\mathcal{B}} 
	\label{eq:optimal_beta_diagonal}\\
	\beta_{ji}(\lambda, \nu) &= &\left[ 1 - \frac{\lambda_i } { \nu_{ij} } \right]_{\mathcal{B}}
	\label{eq:optimal_beta_offdiagonal}
\end{eqnarray}
for all $i \neq j \in \{1,\ldots,m\}$, where $[\; ]_{\mathcal{B}}$ denotes the projection to the set defined in \eqref{eq:beta_set_def}.

\section{Proof of Theorem~\ref{thm:convergence}}

%\begin{pf*}{Proof.}
%
A sufficient condition for \eqref{eq:alg_constraint_conv} and \eqref{eq:alg_value_conv} is that 
\begin{eqnarray}
	&\lim_{t \rightarrow \infty}& \, \mathbb{E}_{h_i} [\alpha_{i}(h_i;t)] = \mathbb{E}_{h_i} [\alpha_{i}^*(h_i)] 
	\label{eq:alpha_convergence_1}\\
	&\lim_{t \rightarrow \infty}& \, \mathbb{E}_{h_i} [\alpha_{i}(h_i;t) q(h_i) ]= \mathbb{E}_{h_i} [\alpha_{i}^*(h_i) q(h_i)]
	\label{eq:alpha_convergence_2}
\end{eqnarray}
hold for all $i =1,\ldots,m$. Indeed this immediately implies \eqref{eq:alg_value_conv}, while \eqref{eq:alg_constraint_conv} is also implied since the optimal policy $\alpha^*$ for problem \eqref{eq:main_optimization}-\eqref{eq:main_constraint} needs to be feasible. In the proof of Theorem~\ref{thm:optimal_characterization} we argued that problem \eqref{eq:main_optimization}-\eqref{eq:main_constraint} is equivalent to the auxiliary problem \eqref{eq:primal_optimization}-\eqref{eq:third_ineq} with variables $\alpha \in \mathcal{A}$, $\beta \in \mathcal{B}^{m \times m}$. Hence it suffices to show that the algorithm converges to the optimal solution of this auxiliary problem in the sense of \eqref{eq:alpha_convergence_1}-\eqref{eq:alpha_convergence_1}.

Recall that after introducing dual variables $\lambda \in \reals_+^m$, $\nu \in \reals_+^{m \times m}$, the Lagrange dual function $g(\lambda, \nu)$ of the auxiliary problem \eqref{eq:primal_optimization}-\eqref{eq:third_ineq} is defined in \eqref{eq:dual_function}. We begin by arguing that that at each iteration of the algorithm, the dual variables $\lambda(t)$, $\nu(t)$ according to \eqref{eq:dual_step_nu_diag}-\eqref{eq:dual_step_lambda} move towards a subgradient direction of the dual function. For convenience let us denote the direction of the steps at \eqref{eq:dual_step_nu_diag}-\eqref{eq:dual_step_nu_offdiag} by the matrix $s_\nu(t) \in \reals^{m \times m}$ defined as
\begin{eqnarray}\label{eq:diagonal_subgradient}
	s_{\nu, ii}(t) &=& \beta_{ii}(t) - \mathbb{E}_{h_i} [ \alpha_i(h_i; t) \, q(h_i)] \\
	\label{eq:offdiagonal_subgradient}
	s_{\nu, ij}(t) &=& \mathbb{E}_{h_j} [\alpha_j(h_j; t)] q_{ji} - \beta_{ji}(t)
\end{eqnarray}
for all $i \neq j \in \{1, \ldots, m\}$, and the steps at \eqref{eq:dual_step_lambda} by the vector $s_\lambda(t) \in \reals^m$ defined as
\begin{eqnarray}\label{eq:vector_subgradient}
	s_{\lambda, i}(t) = \log(c_i) -  \log(  \beta_{ii}(t)) %\, \notag\\
	- \sum_{j\neq i} \log ( 1 - \beta_{ji}(t) ) 
\end{eqnarray}
for all $i \in \{1, \ldots, m\}$. We argue that $s_\nu(t), s_\lambda(t)$ are subgradient directions for the dual function at the point $\lambda(t), \nu(t)$, i.e., that 
\begin{eqnarray}\label{eq:subgradient_defn}
	g(\lambda', \nu') - g(\lambda(t), \nu(t)) \leq \,
	&&(\,\lambda' - \lambda(t)\,)^T s_\lambda(t) \notag\\
	&&+ \, \text{Tr}(\, (\nu'-\nu(t)) \, s_\nu(t) \,)
\end{eqnarray}
for all $\lambda' \in \reals_+^m$, $\nu' \in \reals_+^{m \times m}$. This can be shown as follows. 

%****The next step of the algorithm is to compute the matrix $s_\nu(t) \in \reals^{m \times m}$ and the vector $s_\lambda(t) \in \reals^m$ in \eqref{eq:diagonal_subgradient}-\eqref{eq:vector_subgradient}. We observe that these expressions correspond to the feasibility gap of the variables $\alpha(t)$, $\beta(t)$ computed before with respect to the constraints \eqref{eq:first_ineq}-\eqref{eq:third_ineq} of the auxiliary problem. 

Consider an iteration of Algorithm~\ref{alg:random_access}. The variable $\alpha(t)$ selected by the algorithm at step \eqref{eq:primal_step_alpha} is a variable that minimizes the Lagrangian $L(\alpha, \beta, \lambda(t), \nu(t))$ with respect to the variable $\alpha \in \mathcal{A}$. This follows directly from Lemma~\ref{lemma:Lagrange_minimizers}. Similarly the variables $\beta(t)$ at step \eqref{eq:primal_step_beta_offdiagonal} minimize the Lagrangian function $L(\alpha, \beta, \lambda(t), \nu(t)$ with respect to the variable $\beta \in \mathcal{B}^{m \times m}$. This fact is included in the proof of Lemma~\ref{lemma:Lagrange_minimizers} at \eqref{eq:optimal_beta_diagonal}-\eqref{eq:optimal_beta_offdiagonal}. As a result by the definition of the dual function in \eqref{eq:dual_function} it follows that $g(\lambda(t), \nu(t)) = L(\alpha(t), \beta(t), \lambda(t), \nu(t))$. Additionally we can substitute the Lagrangian at the right hand side with the form given at \eqref{eq:Lagrangian} to get
\begin{eqnarray}\label{eq:dual_fn_1}
	g(\lambda(t), \nu(t)) = \sum_{i=1}^m \; \mathbb{E}_{h_i} \alpha_{i}(h_i;t) p_{i} 
	+ \lambda(t)^T s_\lambda(t) \notag\\
	+ \text{Tr}( \nu(t) s_\nu(t)) .
\end{eqnarray}
Here for convenience we replaced the lengthy parentheses of \eqref{eq:Lagrangian} by the equivalent terms $s_\lambda(t)$, $s_\nu(t)$ defined in  \eqref{eq:offdiagonal_subgradient}, \eqref{eq:vector_subgradient}, \eqref{eq:diagonal_subgradient}.

Next note that at any point $\lambda', \nu'$ the dual function $g(\lambda', \nu')$ is by definition \eqref{eq:dual_function} the minimum of the Lagrangian $L(\alpha, \beta, \lambda', \nu')$, hence we must have that $ g(\lambda', \nu') \leq L(\alpha(t), \beta(t), \lambda', \nu')$. Using again the notation $s_\lambda(t), s_\nu(t)$ at the right hand side we get that
\begin{eqnarray}\label{eq:dual_fn_2}
	g(\lambda', \nu') \leq \sum_{i=1}^m \; \mathbb{E}_{h_i} \alpha_{i}(h_i;t) p_{i}  
	+ \lambda'^T s_\lambda(t) + \text{Tr}( \nu' s_\nu(t)),
	\notag \\
\end{eqnarray}
Subtracting \eqref{eq:dual_fn_1} from \eqref{eq:dual_fn_2} by sides yields \eqref{eq:subgradient_defn}.

To sum up, at each iteration of the algorithm the dual variables $\lambda(t)$, $\nu(t)$ move towards a subgradient direction of the dual function. Additionally the subgradients are bounded. That is true for $s_\nu(t)$ because all the terms at the right hand side of \eqref{eq:diagonal_subgradient}-\eqref{eq:offdiagonal_subgradient} are between 0 and 1. It is also true for $s_\lambda(t)$ because the logarithms at the right hand side of \eqref{eq:vector_subgradient} are finite by the restriction $\beta(t) \in \mathcal{B}^{m \times m}$ defined in \eqref{eq:beta_set_def}. Under the bounded subgradient condition, convergence of $\lambda(t)$, $\nu(t)$ to the optimal dual variables $\lambda^*$, $\nu^*$ for stepsizes as in the statement of the theorem relies on standard subgradient method arguments -- see, e.g., \cite[Prop. 8.2.6]{Bertsekas_convex} for a proof. 

In the rest of the proof, based on the established convergence of the dual variables to the optimal ones, we will show that the same holds for the primal variable $\alpha(.;t)$ in the sense of \eqref{eq:alpha_convergence_1}-\eqref{eq:alpha_convergence_1}. Note that at any iteration $t$ the function $\alpha_i(.;t)$ takes the value $1$ when $\nu_{ii}(t) \, q(h_i) \geq p_{i}  + \sum_{j\neq i} \nu_{ji}(t) q_{ij}$ and $0$ otherwise. Due to the strict monotonicity of the function $q(.)$ this is a threshold-like function taking the value $1$ when $h_i \geq \bar{h}_i(t) =  q^{-1}( p_{i}  + \sum_{j\neq i} \nu_{ji}(t) q_{ij})/ \nu_{ii}(t))$. Since we have established that $\nu(t) \rightarrow \nu^*$, and since the function $q(.)$ is continuous hence its inverse too, we conclude that the threshold $\bar{h}_i(t)$ converges to $\bar{h}_i^* = q^{-1}( p_{i}  + \sum_{j\neq i} \nu_{ji}^* q_{ij})/ \nu_{ii}^*)$. By Theorem~\ref{thm:optimal_characterization} this limit value equals the threshold of the optimal access policy, which is also of a threshold form.

Hence we conclude that $\alpha_i(.;t) \rightarrow \alpha_i^*(.)$ pointwise for all $h_i \in \reals_+$ except perhaps for the point $h_i^*$, i.e., the optimal threshold point. Since the probability measure $\phi_i$ is non-atomic the point $h_i^*$ has a probability measure zero. Hence $\alpha(.;t) \rightarrow \alpha^*(.)$ almost everywhere. Also both sequences of functions $\alpha_i(.;t)$ and $\alpha_i(.;t) q(.)$ are uniformly bounded in $[0,1]$. By the bounded convergence theorem \cite[Theorem 16.5]{Billingsley} we conclude that convergence in expectation, i.e., \eqref{eq:alpha_convergence_1} and \eqref{eq:alpha_convergence_1}, also holds. 
%%
%\begin{eqnarray}
	%\lim_{t \rightarrow \infty} \, \mathbb{E}_{h_i} \alpha_{i}(h_i;t) = \mathbb{E}_{h_i} \alpha_{i}^*(h_i) \\
	%\lim_{t \rightarrow \infty} \, \mathbb{E}_{h_i} [\alpha_{i}(h_i;t) q(h_i) ]= \mathbb{E}_{h_i} [\alpha_{i}^*(h_i) q(h_i)]
%\end{eqnarray}
%%
This completes the proof. 
\qedsymbol
%
%\end{pf*}

\end{document}